\documentclass{article}
\usepackage{amssymb}
\usepackage{amsfonts}
\usepackage{amsmath}

\newtheorem{theorem}{Theorem}

\newtheorem{lemma}[theorem]{Lemma}

\newtheorem{remark}[theorem]{Remark}

\begin{document}

\title{Branching processes in random environment which extinct at a given
moment\thanks{%
This paper is part of a project supported by the German Research Foundation
(DFG) and the Russian Foundation of Basic Research (Grant DFG-RFBR
08-01-91954)}}
\author{\textsc{C. B\"oinghoff} \thanks{%
Fachbereich Mathematik, Universit\"at Frankfurt, Fach 187, D-60054 Frankfurt
am Main, Germany , boeinghoff@math.uni-frankfurt.de} \hspace{.8cm} \textsc{%
E.E. Dyakonova}\thanks{%
Department of Discrete Mathematics, Steklov Mathematical Institute, 8 Gubkin
Street, 119\,991 Moscow, Russia, elena@mi.ras.ru} \\
%EndAName
\textsc{G. Kersting} \thanks{%
Fachbereich Mathematik, Universit\"at Frankfurt, Fach 187, D-60054 Frankfurt
am Main, Germany, kersting@math.uni-frankfurt.de} \hspace{.8cm} \textsc{V.A.
Vatutin}\thanks{%
Department of Discrete Mathematics, Steklov Mathematical Institute, 8 Gubkin
Street, 119\,991 Moscow, Russia, vatutin@mi.ras.ru }}
\date{}
\maketitle

\begin{abstract}
Let $\left\{ Z_{n},n\geq 0\right\} $ be a critical branching process in
random environment and let $T$ be its moment of extinction. Under the
annealed approach we prove, as $n\rightarrow \infty ,$ a limit theorem for
the number of particles in the process at moment $n$ given $T=n+1$ and a
functional limit theorem for the properly scaled process $\left\{
Z_{nt},\delta \leq t\leq 1-\delta \right\} $ given $T=n+1$ and $\delta \in
(0,1/2)$.
\end{abstract}

{\small \emph{Key words and phrases.} Branching process, random environment,
random walk, change of measure, survival probability, functional limit
theorem\newline
\emph{MSC 2000 subject classifications.} Primary 60J80, Secondary 60G50,
60F17 }

\section{Introduction and main results}

The model of branching processes in random environment which we are dealing
with in this paper was introduced by Smith and Wilkinson \cite{SmWil}. To
give a formal definition of these \ processes denote $\mathcal{M}$ the space
of probability measures on $\mathbb{N}_{0}:=\{0,1,2,...\}$ and let $Q$ be a
random variable taking values in $\mathcal{M}$. An infinite sequence $\Pi
=(Q_{1},Q_{2},\ldots )$ of i.i.d. copies of $Q$ is said to form a \emph{%
random environment}. A sequence of $\mathbb{N}_{0}$-valued random variables $%
Z_{0},Z_{1},\ldots $ is called a \emph{branching process in the random
environment} $\Pi $, if $Z_{0}$ is independent of $\Pi $ and, given $\Pi ,$
the process $Z=(Z_{0},Z_{1},\ldots )$ is a Markov chain with 
\begin{equation}
\mathcal{L}\left( Z_{n}\;|\;Z_{n-1}=z_{n-1},\,\Pi =(q_{1},q_{2},\ldots
)\right) \ =\ \mathcal{L}\left( \xi _{n1}+\cdots +\xi _{nz_{n-1}}\right)
\label{transition}
\end{equation}%
for every $n\geq 1,\,z_{n-1}\in \mathbb{N}_{0}$ and $q_{1},q_{2},\ldots \in 
\mathcal{M}$, where $\xi _{n1},\xi _{n2},\ldots $ are i.i.d. random
variables with distribution $q_{n}$. We can write this as 
\begin{equation}
Z_{n}:=\sum_{i=1}^{Z_{n-1}}\xi _{ni},  \label{defZ}
\end{equation}%
where, given the environment, $Z$ is an ordinary inhomogeneous Galton-Watson
process. Thus, $Z_{n}$ is the $n$th generation size of the population and $%
Q_{n}$ is the distribution of the number of children of an individual at
generation $n-1$. We will denote the corresponding probability measure on
the underlying probability space by $\mathbf{P}$.

In what follows \ we identify $Q$ and $Q_{n},n=1,2,...,$ with (random)
generating functions 
\begin{equation*}
f(s)\ :=\ \sum_{i=0}^{\infty }s^{i}\,Q(\{i\})=:\mathbf{E}\left[ s^{\xi }|Q%
\right] ,\qquad 0\leq s\leq 1\ ,
\end{equation*}%
and%
\begin{equation*}
f_{n}(s)\ :=\ \sum_{i=0}^{\infty }s^{i}\,Q_{n}(\{i\})=:\mathbf{E}\left[
s^{\xi _{n}}|Q_{n}\right] ,\qquad 0\leq s\leq 1\ ,
\end{equation*}%
and make no difference between the tuples $\Pi =(Q_{1},Q_{2},\ldots )$ and $%
\mathbf{f}=(f_{1},f_{2},\ldots )$.

Let%
\begin{equation*}
X_{k}:=\ln f_{k}^{\prime }\left( 1\right) ,\qquad \eta _{k}:=\frac{%
f_{k}^{\prime \prime }(1)}{2\left( f_{k}^{\prime }(1)\right) ^{2}},
\end{equation*}%
and%
\begin{equation*}
\left( f,X,\eta \right) \overset{d}{=}\left( f_{1},X_{1},\eta _{1}\right) .
\end{equation*}%
The sequence 
\begin{equation*}
S_{0}:=0,\;S_{n}:=X_{1}+...+X_{n},\,n\in \mathbb{N}:=\{1,2,...\}
\end{equation*}%
is called the \textit{associated random walk} of the corresponding branching
process in random environment (BPRE).

Following \cite{4h} we call a BPRE critical if $\lim \sup_{n\rightarrow
\infty }S_{n}=+\infty $ $\ $\ and $\lim \inf_{n\rightarrow \infty
}S_{n}=-\infty $ both with probability 1.

Let%
\begin{equation*}
T=\min \left\{ k\geq 0:Z(k)=0\right\}
\end{equation*}%
be the extinction moment of the critical BPRE. The aim of the paper is to
study, as $n\rightarrow \infty ,$ the behavior of the process $\left\{
Z(m),1\leq m\leq n\right\} $ given $T=n+1$. Critical BPRE's conditioned on
extinction at a given moment were investigated in \cite{VD} and \cite{VW09}
under the annealed approach and in \cite{VK08} under the quenched approach.
In all the papers it is assumed that the functions $f_{n}\left( s\right) $
are fractional-linear, namely, 
\begin{equation*}
\frac{1}{1-f_{n}\left( s\right) }=\frac{e^{-X_{n}}}{1-s}+\eta _{n}.
\end{equation*}

In \cite{VW09} the asymptotic behavior of the probability $\mathbf{P}\left(
T=n\right) $ as $n\rightarrow \infty $ is found and a conditional functional
limit theorem for the properly scaled process $\left\{ Z(m),1\leq m\leq
n\right\} $ given $T=n+1$ is proved under the assumption that the
distribution of $X$ belongs to the domain of attraction of a stable law with
parameter $\alpha \in (0,2).$ It was shown that in this case the phenomena
of sudden extinction of the process takes place. Namely, if the process
survives for a long time $(T=n+1\rightarrow \infty )$ then $\log Z_{\left[ nt%
\right] }$ grows, roughly speaking, as $n^{1/\alpha }$ up to moment $n$ and
then the process instantly dies out. In particular, $\log Z_{n}$ is of order 
$n^{1/\alpha }$. This may be interpreted as the evolution of the process in
a favorable environment up to moment $n$ and sudden extinction of the
population at moment $T=n+1\rightarrow \infty $ because of a very
unfavorable, even "catastrophic" environment at moment $n$.

For the case $\mathbf{E}X^{2}<\infty $ the asymptotic behavior of the
probability $\mathbf{P}\left( T=n+1\right) $ as $n\rightarrow \infty $ was
investigated in \cite{VD}$.$ However, no functional limit theorem was
proved. We fill this gap \ in the present paper and establish a conditional
functional limit theorem for the process 
\begin{equation*}
\left\{ Z_{nt}e^{-S_{nt}},\delta \leq t\leq 1-\delta \,\big|\,T=n+1\right\}
,\,\delta \in \left( 0,1/2\right) ,
\end{equation*}%
as $n\rightarrow \infty $ and, in addition, show that the conditional law $%
\mathcal{L}\left( Z_{n}\,|\,T=n+1\right) $ weakly converges to a law
concentrated on natural numbers. Thus, contrary to the case considered in 
\cite{VW09}, the phenomenon of sudden extinction is absent if $\mathbf{E}%
X^{2}<\infty $ under the annealed approach.

Note that paper \cite{VK08} demonstrates that in case of the quenched
approach the phenomenon of sudden extinction does not occur if $X$ belongs
to the domain of attraction of a stable law with parameter $\alpha \in (0,2]$%
.

Now we list the basic conditions imposed in this paper on the
characteristics of our BPRE.

\textbf{Assumption A1. }There exists a constant $\chi \in (0,1/2)$ such that 
\textbf{\ }%
\begin{equation}
0<\chi \leq f(0)\leq 1-\chi ,\quad \eta \geq \chi  \label{Asss1}
\end{equation}%
with probability 1.

\textbf{Assumption A2.} The distribution of $X$ has zero mean, finite and
positive variance $\sigma ^{2}$ and is non-lattice.

Let%
\begin{equation*}
\zeta \left( a\right) :=\sum_{y=a}^{\infty }y^{2}Q\left( \left\{ y\right\}
\right) /m\left( Q\right) ^{2},\,a\in \mathbb{N}:=\left\{ 1,2,...\right\} .
\end{equation*}

\textbf{Assumption A3}. For some $\varepsilon >0$ and some $a\in \mathbb{N}$%
\begin{equation*}
\mathbf{E}\left[ \left( \log ^{+}\zeta \left( a\right) \right)
^{2+\varepsilon }\right] <\infty ,
\end{equation*}%
where $\log ^{+}x:=\log \left( \max \left\{ x,1\right\} \right) .$

Here are our main results.

\begin{theorem}
\label{disc} Under A1 to A3, as $n\rightarrow \infty ,$%
\begin{equation*}
\mathbf{P}\left( T=n+1\right) \thicksim cn^{-3/2},
\end{equation*}%
where $c\in (0,\infty ),$ and%
\begin{equation*}
\mathcal{L}\left( Z_{n}\,|\,T=n+1\right) \rightarrow \mathcal{L}\left(
Y\right)
\end{equation*}%
weakly, where $Y$ is a non-degenerate random variable finite with
probability 1. \ 
\end{theorem}

\begin{theorem}
\bigskip \label{fucttheor} Under A1 to A3 for any $\delta \in \left(
0,1/2\right) ,$ as$\ n\rightarrow \infty ,$%
\begin{equation*}
\mathcal{L}\left( \frac{Z_{nt}}{e^{S_{nt}}},t\in \left[ \delta ,1-\delta %
\right] \,\Big|\,T=n+1\right) \Rightarrow \mathcal{L}\left( W_{t},t\in \left[
\delta ,1-\delta \right] \right) ,
\end{equation*}%
where the limiting process $W_{t}$ has a.s. constant trajectories, i.e.,%
\begin{equation*}
\mathbf{P}\left( W_{t}=W\text{ for all }t\in (0,1)\right) =1,
\end{equation*}%
and%
\begin{equation*}
\mathbf{P}\left( 0<W<\infty \right) =1.
\end{equation*}
\end{theorem}

\bigskip Here the symbol $\Rightarrow $ means weak convergence with respect
to the Skorokhod topology in the space $D[\delta ,1-\delta ]$ of cadlag
functions on the interval $\left[ \delta ,1-\delta \right] $.

The proofs of the above results are based on the approach initiated in \cite%
{4h} and developed recently in \cite{ABKV} and use the fact that the
asymptotic behavior of the critical BPRE's is, essentially, specified by the
properties of its associated random walk.

\section{Some auxiliary results}

\bigskip In this section we give a list of general results related with an
oscillating random walk $S_{0},S_{k}=S_{0}+X_{1}+...+X_{k}$ with no
referring to the critical BPRE's and we allow here $S_{0}$ to be a random
variable for technical reason. These results are basically taken from \cite%
{ABKV} and are established under the following assumption.

\textbf{Assumption A4}. There are numbers $c_{n}\rightarrow \infty $ such
that the sequence $S_{n}/c_{n}$ converges in distribution to an $\alpha -$
stable law which is neither concentrated on $\mathbb{R}_{+}:=[0,\infty )$
and $\mathbb{R}_{-}:=(-\infty ,0]$. It is nonlattice.

Introduce the random variables 
\begin{equation*}
M_{n}:=\max \left( S_{1},...,S_{n}\right) ,\quad L_{n}:=\min \left(
S_{1},...,S_{n}\right)
\end{equation*}%
and, given $S_{0}=0,$ the right-continuous functions $u:\mathbb{R}$ $%
\rightarrow \mathbb{R}_{+}$ and $v:$ $\mathbb{R}$ $\rightarrow \mathbb{R}%
_{+} $ specified by the equalities%
\begin{equation*}
u(x):=1+\sum_{k=1}^{\infty }\mathbf{P}\left( -S_{k}\leq x,M_{k}<0\right)
,\quad x\in \mathbb{R}_{+},
\end{equation*}%
\begin{equation*}
v(x):=1+\sum_{k=1}^{\infty }\mathbf{P}\left( -S_{k}>x,L_{k}\geq 0\right)
,\quad x\in \mathbb{R}_{-},
\end{equation*}%
with $u(0)=v(0)=1$ and $u(-x)=v(x)=0,\quad x\in \mathbb{R}_{+}.$

One may check (see, for instance, \cite{4h}\ and \cite{ABKV} ) that for any
oscillating random walk 
\begin{equation}
\mathbf{E}\left[ u(x+X);X+x\geq 0\right] =u(x),\quad x\in \mathbb{R}_{+},
\label{Mes1}
\end{equation}%
\begin{equation}
\mathbf{E}\left[ v(x+X);X+x<0\right] =v(x),\quad x\in \mathbb{R}_{-}.
\label{Mes2}
\end{equation}%
By $u$ and $v$ we construct two probability measures $\mathbf{P}^{+}$ and $%
\mathbf{P}^{-}$. To this aim let $O_{1},O_{2},...$ be a sequence of
identically distributed random variables in a state space $\mathcal{D}$,
adapted to a filtration $\left( \mathcal{F}_{n},n\in \mathbb{N}_{0}\right) $
(possibly larger than the filtration generated by $\left( O_{n},n\geq
1\right) $) such that for all $n$, $O_{n+1}$ is independent of $\mathcal{F}$
and, in particular, $\left( O_{n},n\geq 1\right) $ is a sequence of i.i.d.
random variables. Let, further, $R_{0},R_{1},...$ \ be a sequence of random
variables in a state space $\mathcal{S}$ and also adapted to $\mathcal{F}$.
We assume that the increments $\left( X_{n},n\geq 1\right) $ of the random
walk $S$ are such that for all $n,$ $X_{n}$ are measurable with respect to
the $\sigma $-field generated by $O_{n}$ and $S_{0}$ is $\mathcal{F}_{0}$%
-measurable.

Now for any bounded and measurable function $g:\mathcal{S}\rightarrow 
\mathbb{R}$, we construct probability measures $\mathbf{P}_{x}^{+},x\geq 0,$
and $\mathbf{P}_{x}^{-},x\leq 0,$ fulfilling for each $n$ the equalities%
\begin{equation*}
\mathbf{E}_{x}^{+}\left[ g\left( R_{0},...,R_{n}\right) \right] =\frac{1}{%
u(x)}\mathbf{E}_{x}\left[ g\left( R_{0},...,R_{n}\right) u(S_{n});L_{n}\geq 0%
\right]
\end{equation*}%
and%
\begin{equation*}
\mathbf{E}_{x}^{-}\left[ g\left( R_{0},...,R_{n}\right) \right] =\frac{1}{%
v(x)}\mathbf{E}_{x}\left[ g\left( R_{0},...,R_{n}\right) v(S_{n});M_{n}<0%
\right] .
\end{equation*}

Using (\ref{Mes1})-(\ref{Mes2}) it is not difficult to check (see \cite{4h}
and \cite{ABKV} for more detail) that the measures specified in this way are
consistent in $n$.~

Let $d_{n}=\left( nc_{n}\right) ^{-1}$. In the sequel if no otherwise is
stated, \ we write $a_{n}\sim b_{n}$ if $\lim_{n\rightarrow \infty
}a_{n}/b_{n}=1$ and $a_{n}\rightarrow a$ if $\lim_{n\rightarrow \infty
}a_{n}=a$.

Let 
\begin{equation*}
\tau (n)=\min \left\{ i\leq n:S_{i}=\min \left( S_{0},S_{1},...,S_{n}\right)
\right\}
\end{equation*}%
be the moment of the first random walk minimum up to time $n$.

The next three results are borrowed from \cite{ABKV}.

\begin{lemma}
\label{LasS}( \cite{ABKV}, Proposition 2.1) Under assumption \textbf{A4} for 
$x\geq 0,$ $\theta >0$%
\begin{equation*}
\mathbf{E}_{x}\left[ e^{-\theta S_{n}};L_{n}\geq 0\right] \sim
s(0)d_{n}u(x)\int_{0}^{\infty }e^{-\theta z}v(-z)dz
\end{equation*}%
and for $x\leq 0$%
\begin{equation*}
\mathbf{E}_{x}\left[ e^{\theta S_{n}};\tau (n)=n\right] \sim
s(0)d_{n}v(x)\int_{0}^{\infty }e^{-\theta z}u(z)dz.
\end{equation*}%
In particular, if $\sigma ^{2}:=\mathbf{E}X^{2}<\infty $ then, for some
positive constants $K_{1}$ and $K_{2}$%
\begin{equation}
\mathbf{E}\left[ e^{-S_{n}};L_{n}\geq 0\right] \sim K_{1}n^{-3/2}\text{ and
\thinspace }\mathbf{E}\left[ e^{S_{n}};\tau (n)=n\right] \sim K_{2}n^{-3/2}%
\text{.}  \label{Star}
\end{equation}
\end{lemma}

For $\theta >0,$ let $\mu _{\theta },\nu _{\theta }$ be the probability
measures on $\mathbb{R}_{+}$ and $\mathbb{R}_{-}$ given by their densities 
\begin{equation}
\mu _{\theta }\left( dz\right) :=c_{1\theta }e^{-\theta z}u(z)I\left( z\geq
0\right) dz,\quad \nu _{\theta }\left( dz\right) :=c_{2\theta }e^{\theta
z}\nu (z)I\left( z<0\right) dz,  \label{DefMeasure}
\end{equation}%
where%
\begin{equation}
c_{1\theta }^{-1}=\int_{0}^{\infty }e^{-\theta z}u(z)dz,\quad c_{2\theta
}^{-1}=\int_{0}^{\infty }e^{-\theta z}v(-z)dz.  \label{DefConst}
\end{equation}

\begin{lemma}
\label{Lminimal}( \cite{ABKV}, Proposition 2.7) Let $\ 0<\delta <1.$ Let $%
U_{n}=g_{n}\left( R_{0},...,R_{[\delta n]}\right) ,$\ $n\geq 1,$ be random
variables with values in an Euclidean (or polish) space $S$ such that, as $%
n\rightarrow \infty $%
\begin{equation*}
U_{n}\rightarrow U_{\infty }\text{ \quad }\mathbf{P}^{+}\text{ a.s.}
\end{equation*}%
for some $S$-valued random variable $U_{\infty }$. Also, let $%
V_{n}=h_{n}\left( Q_{1},...,Q_{[\delta n]}\right) $, $\,n\geq 1,$ be random
variables with values in an Euclidean (or polish) space $S^{\prime }$ such
that%
\begin{equation*}
V_{n}\rightarrow V_{\infty }\text{\quad }\mathbf{P}_{x}^{-}\text{ a.s.}
\end{equation*}%
for all $x\leq 0$ and some $S^{\prime }$-valued random variable $V_{\infty }$%
. Denote%
\begin{equation*}
\tilde{V}_{n}:=h_{n}\left( Q_{n},...,Q_{n-\left[ \delta n\right] +1}\right) .
\end{equation*}%
Under assumption \textbf{A4 }for $\theta >0$ and any bounded continuous
function $\varphi :S\times S^{\prime }\times \mathbb{R\rightarrow R}$ as $%
n\rightarrow \infty $%
\begin{eqnarray*}
\frac{\mathbf{E}\left[ \varphi \left( U_{n},\tilde{V}_{n},S_{n}\right)
e^{-\theta S_{n}};L_{n}\geq 0\right] }{\mathbf{E}\left[ e^{-\theta
S_{n}};L_{n}\geq 0\right] }&& \\
\rightarrow \idotsint &\varphi \left( u,v,-z\right)& \mathbf{P}^{+}\left(
U_{\infty }\in du\right) \mathbf{P}_{z}^{-}\left(V_{\infty }\in dv\right)\nu
_{\theta }\left( dz\right) .
\end{eqnarray*}
\end{lemma}

The dual version of Lemma \ref{Lminimal} looks as follows.

\begin{lemma}
\label{Lend} ( \cite{ABKV}, Proposition 2.9) Let $U_{n},V_{n},\tilde{V}%
_{n},n=1,2,...,\infty $, be as in Lemma \ref{Lminimal} and now fulfilling,
as $n\rightarrow \infty $%
\begin{equation*}
U_{n}\rightarrow U_{\infty }\quad \text{ }\mathbf{P}_{x}^{+}-\text{ a.s., \ }%
V_{n}\rightarrow V_{\infty }\text{ \ \ \ }\mathbf{P}^{-}-\text{ a.s.}
\end{equation*}%
for all $x\geq 0$. Under assumption \textbf{A4} for any bounded continuous
function $\varphi :S\times S^{\prime }\times \mathbb{R\rightarrow R}$ and
for $\theta >0$ as $n\rightarrow \infty $%
\begin{eqnarray*}
&&\frac{\mathbf{E}\left[ \varphi \left( U_{n,}\tilde{V}_{n},S_{n}\right)
e^{\theta S_{n}};\tau (n)=n\right] }{\mathbf{E}\left[ e^{\theta S_{n}};\tau
(n)=n\right] } \\
&&\qquad \qquad \rightarrow \int ...\int \varphi \left( u,v,-z\right) 
\mathbf{P}_{z}^{+}\left( U_{\infty }\in du\right) \mathbf{P}^{-}\left(
V_{\infty }\in dv\right) \mu _{\theta }\left( dz\right) .
\end{eqnarray*}
\end{lemma}

\begin{remark}
\label{R1}It is easy to see (by introducing formal arguments) that the
statements of the lemmas are valid for any integer-valued function $w(n)$
such that $w(n)\leq \delta n$ for all sufficiently large $n$. Then the
functions $g_{n}$ and $h_{n}$ can be viewed as functions also of the missing
variables). Later on we use this fact with no additional reference.
\end{remark}

\section{\protect\bigskip\ Discrete limit distribution\label{SecDisc}}

Introduce the compositions 
\begin{equation}
f_{k,n}(s)\ :=\ f_{k+1}(f_{k+2}(\cdots f_{n}(s)\cdots ))\ ,\ 0\leq k<n,\quad
f_{n,n}(s)\ :=s,  \label{erz}
\end{equation}%
and%
\begin{equation*}
f_{k,0}(s)\ :=\ f_{k}(f_{k-1}(\cdots f_{1}(s)\cdots )),\,k>0.
\end{equation*}%
In this notation we may rewrite the distributional identity (\ref{transition}%
) for $k\leq n$ as 
\begin{equation}
\mathbf{E}[s^{Z_{n}}\;|\;\mathbf{f},Z_{k}]=\mathbf{E}[s^{Z_{n}}\;|\;\Pi
,Z_{k}]\ =\ f_{k,n}(s)^{Z_{k}}\quad \mathbf{P}\text{--a.s.}
\label{branching2}
\end{equation}

For $0\leq k\leq n$ and $S_{0}:=0$ let%
\begin{equation*}
a_{k,n}:=e^{-(S_{n}-S_{k})},\quad a_{n}:=a_{0,n}=e^{-S_{n}},
\end{equation*}%
\begin{equation}
b_{k,n}:=\sum_{i=k}^{n-1}\eta _{i+1}e^{-(S_{i}-S_{k})},\quad
b_{n}:=b_{0,n}=\sum_{i=0}^{n-1}\eta _{i+1}e^{-S_{i}}.  \label{DefBn}
\end{equation}

\begin{lemma}
\label{LfrLin}(see, for instance, \cite{GK}) In the fractional-linear case
for any $0\leq j<n$%
\begin{equation}
\left( 1-f_{j,n}\left( s\right) \right) ^{-1}=\frac{a_{j,n}}{1-s}+b_{j,n}.
\label{base}
\end{equation}%
In particular,%
\begin{equation}
\left( 1-f_{j,n}\left( 0\right) \right) ^{-1}=a_{j,n}+b_{j,n}  \label{sim1}
\end{equation}%
and%
\begin{eqnarray}
\left( 1-f_{0,n}\left( s\right) \right) ^{-1} &=&\frac{a_{n}}{1-s}%
+b_{n}=\left( 1-f_{0,j}\left( f_{j,n}\left( s\right) \right) \right) ^{-1} 
\notag \\
&=&\frac{a_{j}}{1-f_{j,n}\left( s\right) }+b_{j}=\frac{a_{j}a_{j,n}}{1-s}%
+b_{j}+a_{j}b_{j,n}.  \label{sim2}
\end{eqnarray}
\end{lemma}

\begin{lemma}
\bigskip \label{Lold1}(see Lemma 3.1 in \cite{ABKV}) If conditions A3 - A4
are valid then for any $x\geq 0$%
\begin{equation*}
B^{+}:=\lim_{n\rightarrow \infty }b_{n}=\sum_{i=0}^{\infty }\eta
_{i+1}e^{-S_{i}}<\infty \quad \mathbf{P}_{x}^{+}\text{- a.s.}
\end{equation*}%
and for any $x\leq 0$%
\begin{equation*}
B^{-}:=\lim_{n\rightarrow \infty }\sum_{i=1}^{n}\eta
_{i}e^{S_{i}}=\sum_{i=1}^{\infty }\eta _{i}e^{S_{i}}<\infty \quad \mathbf{P}%
_{x}^{-}\text{- a.s.}
\end{equation*}
\end{lemma}

Let $\kappa :\mathbb{R}\rightarrow \mathbb{R}_{+}$ be the function specified
by the equality%
\begin{equation*}
\kappa \left( y\right) :=yI\left( y>0\right) ,
\end{equation*}%
where $I\left( A\right) $ is the indicator of the event $A,$ and let, for
positive constants $\alpha ,\beta ,\gamma $ 
\begin{equation}
\phi \left( \alpha ,\beta ,\gamma ;u,v,x\right) :=\frac{1}{e^{-x}\alpha
+\beta +\gamma (\kappa \left( u\right) +e^{-x}\kappa \left( v\right) )},
\label{DefphiDisc}
\end{equation}%
and%
\begin{equation*}
\Phi \left( u,v,x\right) =\Phi \left( \alpha _{1},\alpha _{2},\beta
_{1},\beta _{2},\gamma _{1},\gamma _{2};u,v,x\right) :=\prod_{i=1}^{2}\phi
\left( \alpha _{i},\beta _{i},\gamma _{i};u,v,x\right) .
\end{equation*}%
Clearly, $\Phi \left( u,v,x\right) $ is continuous in $\mathbb{R}^{3}$ and
bounded by $\beta _{1}^{-1}\beta _{2}^{-1}$.

\begin{lemma}
\label{Ldiscleft}Under the conditions A3-A4, for any $\alpha _{i}>0,\beta
_{i}>0,\gamma _{i}>0,i=1,2$ 
\begin{eqnarray*}
&&\mathbf{E}\left[ \frac{e^{-S_{n}}}{\prod_{i=1}^{2}\left( e^{-S_{n}}\alpha
_{i}+\beta _{i}+\gamma _{i}b_{n}\right) };L_{n}\geq 0\right] /\mathbf{E}%
\left[ e^{-S_{n}};L_{n}\geq 0\right] \\
&&\qquad \qquad \quad \rightarrow \iiint \Phi \left( u,v,-z\right) \mathbf{P}%
^{+}\left( B^{+}\in du\right) \mathbf{P}_{z}^{-}\left( B^{-}\in dv\right)
\nu _{1}\left( dz\right) .
\end{eqnarray*}
\end{lemma}

\textbf{Proof}. We have%
\begin{eqnarray*}
e^{-S_{n}}\alpha +\beta +\gamma b_{n} &=&e^{-S_{n}}\alpha +\beta +\gamma
(U_{ \left[ n/2\right] }+e^{-S_{n}}\tilde{V}_{\left[ n/2\right] }) \\
&=&1/\phi \left( \alpha ,\beta ,\gamma ;U_{\left[ n/2\right] },\tilde{V}_{%
\left[ n/2\right] },S_{n}\right) ,
\end{eqnarray*}%
where%
\begin{equation*}
U_{\left[ n/2\right] }=g_{\left[ n/2\right] }\left( Q_{1},...,Q_{\left[ n/2%
\right] }\right) :=b_{\left[ n/2\right] }
\end{equation*}%
with $b_{n}$ specified by (\ref{DefBn}), and 
\begin{equation*}
\tilde{V}_{\left[ n/2\right] }=h_{\left[ n/2\right] }\left( Q_{n},...,Q_{n-%
\left[ n/2\right] +1}\right) :=\sum_{i=0}^{n-\left[ n/2\right] }\eta _{i+%
\left[ n/2\right] +1}e^{S_{n}-S_{i+\left[ n/2\right] }}.
\end{equation*}%
By Lemma \ref{Lold1} as $n\rightarrow \infty $ for any $x\geq 0$%
\begin{equation}
U_{\left[ n/2\right] }=b_{\left[ n/2\right] }\rightarrow B^{+}\quad \mathbf{P%
}_{x}^{+}\text{-a.s.}  \label{Conr}
\end{equation}%
and, for any $x\leq 0$ 
\begin{equation}
V_{\left[ n/2\right] }=\sum_{i=1}^{\left[ n/2\right] +1}\eta
_{i}e^{S_{i}}\rightarrow B^{-}\quad \mathbf{P}_{x}^{-}\text{-a.s.}
\label{Conl}
\end{equation}%
Applying Lemma \ref{Lminimal} to $\Phi \left( U_{n},\tilde{V}%
_{n},S_{n}\right) $ completes the proof of the desired statement.

Let, for $\alpha >0$ 
\begin{equation}
\psi \left( \alpha ;u,v,x\right) :=\frac{1}{\alpha +e^{x}\kappa \left(
u\right) +\kappa \left( v\right) }  \label{DefpsiDisc}
\end{equation}%
and\textbf{\ }%
\begin{equation}
\Psi \left( u,v,x\right) =\Psi \left( \alpha _{1},\alpha _{2};u,v,x\right)
:=\prod_{i=1}^{2}\psi \left( \alpha _{i};u,v,x\right) .  \label{MoreDef}
\end{equation}

\begin{lemma}
\label{Ldiscright}Under the conditions A3-A4, for any $\alpha _{1}>0,\alpha
_{2}>0$ 
\begin{eqnarray*}
&&\mathbf{E}\left[ \frac{e^{S_{n}}}{\prod_{i=1}^{2}\left( \alpha
_{i}+e^{S_{n}}b_{n}\right) };\tau (n)=n\right] /\mathbf{E}\left[
e^{S_{n}};\tau (n)=n\right] \\
&&\qquad \qquad \rightarrow \iiint \Psi \left( u,v,-z\right) \mathbf{P}%
_{z}^{+}\left( B^{+}\in du\right) \mathbf{P}^{-}\left( B^{-}\in dv\right)
\mu _{1}\left( dz\right) \mathbf{.}
\end{eqnarray*}
\end{lemma}

\textbf{Proof.} We have%
\begin{equation*}
\alpha +e^{S_{n}}b_{n}=\alpha +e^{S_{n}}U_{\left[ n/2\right] }+\tilde{V}_{%
\left[ n/2\right] }=1/\psi \left( \alpha ;U_{\left[ n/2\right] },\tilde{V}_{%
\left[ n/2\right] },S_{n}\right) ,
\end{equation*}%
where $U_{\left[ n/2\right] }$ and $\tilde{V}_{\left[ n/2\right] }$ are the
same as in Lemma \ref{Ldiscleft}. Now using (\ref{Conr}) and (\ref{Conl})
once again it is not difficult to complete the proof of the lemma.

\textbf{Proof of Theorem \ref{disc}}. For $s\in \lbrack 0,1)$ denote 
\begin{equation*}
X_{f}(s):=\frac{sf(0)}{1-sf(0)},\quad G_{n}(s):=1-f_{0,n}(s)=\left( \frac{%
a_{n}}{1-s}+b_{n}\right) ^{-1}
\end{equation*}%
and let%
\begin{eqnarray}
\Delta _{n}\left( s\right) := &&f_{0,n}(sf\left( 0\right) )-f_{0,n}\left(
0\right)  \notag \\
&=&\left( 1-f_{0,n}\left( 0\right) \right) (1-f_{0,n}\left( sf(0)\right)
)e^{-S_{n}}\frac{sf(0)}{1-sf(0)}  \label{Interm} \\
&=&G_{n}(f(0))G_{n}(sf(0))X_{f}(s)e^{-S_{n}}  \notag \\
&=&\left( a_{n}+b_{n}\right) ^{-1}\left( \frac{a_{n}}{1-sf(0)}+b_{n}\right)
^{-1}X_{f}(s)e^{-S_{n}},  \notag
\end{eqnarray}%
where we have used the explicit form of $f_{0,n}(s)$ and the equality $f%
\overset{d}{=}f_{n+1}$. It is not difficult to check that%
\begin{eqnarray*}
\mathbf{E}\left[ s^{Z_{n}};T=n+1\right] &=&\mathbf{E}\left[
s^{Z_{n}};Z_{n}>0,Z_{n+1}=0\right] \\
&=&\mathbf{E}\left[ \left( sf_{n+1}(0)\right) ^{Z_{n}}I\left( Z_{n}>0\right) %
\right] \\
&=&\mathbf{E}\left[ \left( sf_{n+1}(0)\right) ^{Z_{n}}-I\left(
Z_{n}=0\right) \right] =\mathbf{E}\Delta _{n}\left( s\right) \\
&=&D_{1}(N,n)+D_{2}(N,n)+D_{3}(N,n),
\end{eqnarray*}%
where%
\begin{equation*}
D_{1}(N,n):=\sum_{j=0}^{N}\mathbf{E}\left[ \Delta _{n}\left( s\right) ;\tau
(n)=j\right] ,
\end{equation*}%
\begin{equation*}
D_{2}(N,n):=\sum_{j=N+1}^{n-N-1}\mathbf{E}\left[ \Delta _{n}\left( s\right)
;\tau (n)=j\right] ,
\end{equation*}%
\begin{equation*}
D_{3}(N,n):=\sum_{j=n-N}^{n}\mathbf{E}\left[ \Delta _{n}\left( s\right)
;\tau (n)=j\right] .
\end{equation*}%
By (\ref{Interm}), the evident inequalities 
\begin{equation*}
1-f_{0,n}\left( 0\right) =\min_{1\leq k\leq n}\left( 1-f_{0,k}\left(
0\right) \right) \leq \min_{1\leq k\leq n}e^{S_{k}}=e^{S_{\tau (n)}},
\end{equation*}%
Assumption A1 and the estimates 
\begin{equation}
X_{f}(s)\leq X_{f}(1)\leq (1-\chi )\chi ^{-1}=:\rho  \label{EstX}
\end{equation}%
following from (\ref{Asss1})\ we obtain%
\begin{equation*}
\Delta _{n}\left( s\right) \leq \left( 1-f_{0,n}\left( 0\right) \right)
^{2}e^{-S_{n}}X_{f}(1)\leq \rho \chi e^{2S_{\tau (n)}-S_{n}}\leq \rho
e^{2S_{\tau (n)}-S_{n}}.
\end{equation*}%
Using this estimate, the asymptotic relation (\ref{Star}) and the duality
principle for random walks it is not difficult to show that for any $%
\varepsilon >0$ one can find $N=N(\varepsilon )$ such that for all
sufficiently large $n\geq 2N+1$ 
\begin{eqnarray}
D_{2}(N,n) &\leq &\rho \sum_{j=N+1}^{n-N-1}\mathbf{E}\left[ e^{2S_{\tau
(n)}-S_{n}};\tau (n)=j\right]  \notag \\
&=&\rho \sum_{j=N+1}^{n-N-1}\mathbf{E}\left[ e^{S_{j}};\tau (j)=j\right] 
\mathbf{E}\left[ e^{-S_{n-j}};L_{n-j}\geq 0\right]  \notag \\
&\leq &const\times \sum_{j=N+1}^{n-N-1}\frac{1}{j^{3/2}}\frac{1}{(n-j)^{3/2}}%
\leq \varepsilon n^{-3/2}.  \label{middl}
\end{eqnarray}

Further, \ for fixed $j$ let%
\begin{equation*}
\tilde{V}_{\left[ (n-j)/2\right] }:=\sum_{i=0}^{n-\left[ \left( n-j\right) /2%
\right] }\eta _{i+\left[ \left( n-j\right) /2\right] +1}e^{S_{n}-S_{i+\left[
\left( n-j\right) /2\right] }}.
\end{equation*}%
By (\ref{sim2}) and (\ref{DefphiDisc}) we have for $s\in \lbrack 0,1)$%
\begin{eqnarray}
G_{n}(sf(0)) &=&\left( \frac{a_{j}a_{j,n}}{1-sf(0)}+b_{j}+a_{j}b_{j,n}%
\right) ^{-1}  \notag \\
&=&\left( \frac{a_{j}}{1-sf(0)}e^{-(S_{n}-S_{j})}+b_{j}+a_{j}\left( b_{j,%
\left[ \left( n-j\right) /2\right] }+e^{-(S_{n}-S_{j})}\tilde{V}_{\left[
(n-j)/2\right] }\right) \right) ^{-1}  \notag \\
&=&\phi \left( \frac{a_{j}}{1-sf(0)},b_{j},a_{j};b_{j,\left[ \left(
n-j\right) /2\right] },\tilde{V}_{\left[ (n-j)/2\right] },S_{n}-S_{j}\right)
\leq 1.  \label{Gestimate}
\end{eqnarray}%
Denote by $\mathcal{F}_{j}^{\ast }$ the $\sigma $-algebra generated by $%
Q,Q_{1},..,Q_{j}$ and introduce a temporary notation%
\begin{equation*}
\alpha _{1}=a_{j},\quad \alpha _{2}=\frac{a_{j}}{1-sf(0)}.
\end{equation*}%
Then%
\begin{equation}
\mathbf{E}\left[ \Delta _{n}\left( s\right) ;\tau (n)=j|\mathcal{F}%
_{j}^{\ast }\right] =\mathbf{E}\left[ a_{j}X_{f}(s)I\left( \tau (j)=j\right)
A_{j,n}(s)\right]  \label{AA1}
\end{equation}%
with 
\begin{eqnarray*}
A_{j,n}(s) &:=&\mathbf{E}\left[ \frac{e^{-\hat{S}_{n-j}}}{%
\prod_{i=1}^{2}\left( \alpha _{i}\hat{a}_{n-j}+b_{j}+a_{j}\hat{b}%
_{n-j}\right) };\hat{L}_{n-j}\geq 0|\mathcal{F}_{j}^{\ast }\right] \\
&=&\mathbf{E}\left[ e^{-\hat{S}_{n-j}}\prod_{i=1}^{2}\phi \left( \alpha
_{i},b_{j},a_{j};\hat{b}_{[n-j/2]},\hat{V}_{\left[ (n-j)/2\right] },\hat{S}%
_{n-j}\right) ;\hat{L}_{n-j}\geq 0|\mathcal{F}_{j}^{\ast }\right]
\end{eqnarray*}%
and where we have taken the agreement that a random variable $\hat{\zeta}=%
\hat{\zeta}\left( \hat{Q}_{1},..,\hat{Q}_{n}\right) $ has the same
definition as $\zeta =\zeta \left( Q_{1},...,Q_{n\text{ }}\right) $ but is
generated by a sequence $\hat{Q}_{1},..,\hat{Q}_{n}$ which is independent of 
$\mathcal{F}_{j}^{\ast }$ and has the same distribution as $Q_{1},...,Q_{n%
\text{ }}$. By (\ref{Gestimate}) and asymptotic representation (\ref{Star})
we conclude that for each $j$ there exists a constant $d_{j}$ such that 
\begin{equation*}
n^{3/2}A_{j,n}(s)\leq n^{3/2}\mathbf{E}\left[ e^{-\hat{S}_{n-j}};\hat{L}%
_{n-j}\geq 0\right] \leq d_{j}
\end{equation*}
for all $n$. Now, using Lemma \ref{Ldiscleft} and the dominated convergence
theorem we see that for each fixed $j$%
\begin{eqnarray}
\mathcal{A}_{j} (s)&:=&\lim_{n\rightarrow \infty }n^{3/2}\mathbf{E}\left[
\Delta _{n}\left( s\right) ;\tau (n)=j\right]  \notag \\
&&=\mathbf{E}\left[ a_{j}X_{f}(s)I\left( \tau (j)=j\right)
\lim_{n\rightarrow \infty }n^{3/2}A_{j,n}(s)\right]  \notag \\
&&\qquad\qquad=K_{1}\mathbf{E}\left[ a_{j}I\left( \tau (j)=j\right)
X_{f}(s)A_{j}(s)\right] ,  \label{Acalif}
\end{eqnarray}%
where%
\begin{eqnarray}
A_{j}(s) &:=&\iiint \Phi \left( a_{j},\frac{a_{j}}{1-sf(0)}%
,b_{j},b_{j},a_{j},a_{j};u,v,-z\right) \times  \notag \\
&&\qquad\qquad\qquad\times \mathbf{P}^{+}\left( B^{+}\in du\right) \mathbf{P}%
_{z}^{-}\left( B^{-}\in dv\right) \nu _{1}\left( dz\right).  \label{Asimp}
\end{eqnarray}

To evaluate $\mathbf{E}\left[ \Delta _{n}\left( s\right) ;\tau (n)=n-j\right]
$ for a fixed $j$ let $\mathcal{\hat{F}}_{j}$ be the the $\sigma $-algebra
generated by a sequence of random laws $Q,\hat{Q}_{1},..,\hat{Q}_{j},$ where 
$\hat{Q}_{1},..,\hat{Q}_{j}$ are distributed as $Q_{1},...,Q_{j}$ and are
independent of $Q_{1},...,Q_{n}$. Introduce a temporary notation%
\begin{equation}
\hat{\alpha}_{1}=\frac{1}{1-\hat{f}_{0,j}\left( 0\right) }\geq 1,\quad \hat{%
\alpha}_{2}=\frac{1}{1-\hat{f}_{0,j}\left( sf(0)\right) }\geq 1.
\label{DefAlpha}
\end{equation}%
By (\ref{sim2}) we see that%
\begin{eqnarray*}
&&\mathbf{E}[\Delta _{n}\left( s\right) ;\tau (n)=n-j] \\
&=&\mathbf{E}\left[ G_{n-j}(f_{n-j,n}\left( 0\right)
)G_{n-j}(f_{n-j,n}\left( sf(0)\right)
e^{-S_{n-j}}e^{-(S_{n}-S_{n-j})}X_{f}(s);\tau (n)=n-j\right] \\
&&\qquad =\mathbf{E}\left[ \hat{a}_{j}I\left( \hat{L}_{j}\geq 0\right)
X_{f}(s)B_{j,n}(s)\right] ,
\end{eqnarray*}%
where%
\begin{equation*}
B_{j,n}(s):=\mathbf{E}\left[ \frac{e^{S_{n-j}}}{\prod_{i=1}^{2}\left( \hat{%
\alpha}_{i}+e^{S_{n-j}}b_{n-j}\right) };\tau (n-j)=n-j|\mathcal{\hat{F}}_{j}%
\right] .
\end{equation*}%
In view of (\ref{DefAlpha}) $\prod_{i=1}^{2}\left( \hat{\alpha}%
_{i}+e^{S_{n-j}}b_{n-j}\right) \geq 1$ which, along with (\ref{Star}),
implies that for each~$j$ there exists a constant $d_{j}^{\prime }$ such
that 
\begin{equation*}
n^{3/2}B_{j,n}(s)\leq n^{3/2}\mathbf{E}\left[ e^{S_{n-j}};\tau (n-j)=n-j%
\right] \leq d_{j}^{\prime }
\end{equation*}%
for all $n$. Now Lemma \ref{Ldiscright}, the duality principle for random
walks and the dominated convergence theorem yield for each~$j$%
\begin{eqnarray*}
\mathcal{B}_{j}(s) &:&=\lim_{n\rightarrow \infty }n^{3/2}\mathbf{E}\left[
\Delta _{n}\left( s\right) ;\tau (n)=n-j\right] \\
&=&\mathbf{E}\left[ \hat{a}_{j}I\left( \hat{L}_{j}\geq 0\right)
X_{f}(s)\lim_{n\rightarrow \infty }n^{3/2}B_{j,n}(s)\right] \\
&&\qquad \quad =K_{2}\mathbf{E}\left[ \hat{a}_{j}I\left( \hat{L}_{j}\geq
0\right) X_{f}(s)B_{j}(s)\right] ,
\end{eqnarray*}%
where (recall (\ref{DefpsiDisc}) and (\ref{MoreDef}))%
\begin{equation*}
B_{j}(s):=\iiint \Psi \left( \hat{\alpha}_{1},\hat{\alpha}_{2};u,v,-z\right) 
\mathbf{P}_{z}^{+}\left( B^{+}\in du\right) \mathbf{P}^{-}\left( B^{-}\in
dv\right) \mu _{1}\left( dz\right)
\end{equation*}%
with $\hat{\alpha}_{1}$ and $\hat{\alpha}_{2}$ specified by (\ref{DefAlpha}%
). As a result letting first $n\rightarrow \infty $ and than $N\rightarrow
\infty $ we get%
\begin{equation*}
H(s):=\lim_{n\rightarrow \infty }n^{3/2}\mathbf{E}\left[ s^{Z_{n}};T=n+1%
\right] =\sum_{j=0}^{\infty }\left( \mathcal{A}_{j}(s)+\mathcal{B}%
_{j}(s)\right) .
\end{equation*}%
In particular,%
\begin{equation}
\lim_{n\rightarrow \infty }n^{3/2}\mathbf{P}\left( T=n+1\right) =H(1).
\label{AsymExt}
\end{equation}%
Hence we conclude that%
\begin{equation*}
\lim_{n\rightarrow \infty }\mathbf{E}\left[ s^{Z_{n}}|T=n+1\right] =\frac{%
H(s)}{H(1)}=:\mathbf{E}\left[ s^{Y}\right]
\end{equation*}%
and, by the dominated convergence theorem and continuity of the functions
involved, $\lim_{s\uparrow 1}H(s)=H(1)$ showing that the limiting
distribution has no atom at infinity.

The theorem is proved.

\section{Functional limit theorem\label{SecFunc}}

The proof of Theorem \ref{fucttheor} uses the same ideas as the proof of
Theorem \ref{disc}.

Let $l:\mathbb{R}\rightarrow \left[ 0,1\right] $ be the function specified
by the equality 
\begin{equation*}
l(y):=yI\left( 0\leq y\leq 1\right) +I\left( y>1\right) .
\end{equation*}%
For parameters $\alpha >0,\beta >0,\lambda >0,$ a three-dimensional vector $%
\mathbf{u}=\left( u_{1},u_{2},u_{3}\right) $ and variables $v$ and $x$
introduce the function 
\begin{equation*}
\theta \left( \alpha ,\beta ,\lambda ;\mathbf{u},v,x\right) :=\frac{l(v)}{%
e^{-x}\alpha +\left( \lambda \alpha e^{-x}l(u_{1})+l(u_{2})l\left( v\right)
\right) \left( \beta +\alpha \kappa (u_{3})\right) }.
\end{equation*}

It is not difficult to check that if $u_{2}\geq \varepsilon $ for some $%
\varepsilon >0$ then $\phi \left( \alpha ,\beta ,\lambda ;\mathbf{u}%
,v,x\right) \leq \beta ^{-1}\varepsilon ^{-1}$ and is continuous in $\mathbf{%
u},v$ and $x$ in the mentioned domain. For the particular case $\mathbf{u}%
=\left( 1,1,u\right) $ we use one more notation%
\begin{equation}
\theta ^{\ast }\left( \alpha ,\beta ,\lambda ;u,v,x\right) :=\frac{l(v)}{%
\alpha e^{-x}+(\lambda \alpha e^{-x}+l(v))\left( \beta +\alpha \kappa
(u)\right) }.  \label{explisphi}
\end{equation}%
With the functions above and $\mathbf{v}=\left( v_{1},v_{2}\right) $ we
associate two more functions%
\begin{eqnarray*}
\Theta \left( \alpha ,\beta ,\lambda ;\mathbf{u},\mathbf{v},x\right)
&:=&\prod_{i=1}^{2}\theta \left( \alpha ,\beta ,\lambda ;\mathbf{u}%
,v_{i},x\right) ,\, \\
\Theta ^{\ast }\left( \alpha ,\beta ,\lambda ;u,\mathbf{v},x\right)
&:=&\prod_{i=1}^{2}\theta ^{\ast }\left( \alpha ,\beta ,\lambda
;u,v_{i},x\right) .
\end{eqnarray*}

For fixed $t\in \left( 0,1\right) $ and $r\in \lbrack 0,1)$ introduce a
random vector$~$%
\begin{equation*}
\ \mathbf{U}_{\left[ nt\right] }=\left( U_{\left[ nt\right] }^{(1)},U_{\left[
nt\right] }^{(2)},U_{\left[ nt\right] }^{(3)}\right)
\end{equation*}%
and a random variable $\tilde{V}_{\left[ nt\right] }(r)$ by the equalities%
\begin{equation}
U_{\left[ nt\right] }^{(1)}:=\frac{1-\exp \left\{ -\lambda \alpha e^{-S_{%
\left[ nt\right] }}\right\} }{\lambda \alpha e^{-S_{\left[ nt\right] }}}%
,\,U_{\left[ nt\right] }^{(2)}:=\exp \left\{ -\lambda \alpha e^{-S_{\left[ nt%
\right] }}\right\} ,U_{\left[ nt\right] }^{(3)}:=b_{\left[ nt\right] },
\label{CompU}
\end{equation}%
and%
\begin{equation}
\tilde{V}_{\left[ nt\right] }(r):=\left( 1-f_{\left[ nt\right] ,n}(r)\right)
e^{S_{\left[ nt\right] }-S_{n}}.  \label{CompV}
\end{equation}%
Let, further,%
\begin{equation*}
V_{k}(r):=\left( 1-f_{k,0}(r)\right) e^{-S_{k}}.
\end{equation*}%
It follows from Lemma \ref{Lold1} that 
\begin{equation*}
\lim_{k\rightarrow \infty }S_{k}=\infty \quad \mathbf{P}_{x}^{+}\text{ -
a.s. for any }x\geq 0.
\end{equation*}%
This and Lemma \ref{Lold1} show that, as $n\rightarrow \infty $ 
\begin{equation}
\mathbf{U}_{n}\rightarrow \mathbf{U}_{\infty }:=\left( 1,1,B^{+}\right)
\quad \mathbf{P}_{x}^{+}\text{ - a.s.}  \label{DDU}
\end{equation}%
Besides, if $\mathbf{\tilde{V}}_{\left[ nt\right] }(r_{1},r_{2}):=\left( 
\tilde{V}_{\left[ nt\right] }(r_{1}),\tilde{V}_{\left[ nt\right]
}(r_{2})\right) $ then as $k\rightarrow \infty $ (see, for instance, formula
(3.1) in \cite{VD2}) 
\begin{eqnarray}
\mathbf{V}_{k}(r_{1},r_{2}) &:=&\left( V_{k}(r_{1}),V_{k}(r_{2})\right) 
\notag \\
&&\rightarrow \left( V_{\infty }(r_{1}),V_{\infty }(r_{2})\right) =:\mathbf{V%
}_{\infty }(r_{1},r_{2})\quad \mathbf{P}_{x}^{-}\text{ - a.s.,}  \label{DDV1}
\end{eqnarray}%
where 
\begin{equation}
1/V_{\infty }(r):=B^{-}+\left( 1-r\right) ^{-1}.  \label{DefVinfin}
\end{equation}%
Observe that 
\begin{equation*}
\mathbf{P}^{+}\left( 0<B^{+}<\infty \right) =\mathbf{P}^{-}\left(
0<V_{\infty }\left( r\right) <1\right) =1
\end{equation*}%
by Lemma \ref{Lold1}.

The arguments above combined with Lemma \ref{Lminimal} imply the following
statement.

\begin{lemma}
\label{LLleft}Under conditions A1 to A4%
\begin{eqnarray*}
&&\frac{\mathbf{E}\left[ \Theta \left( \alpha ,\beta ,\lambda ;\mathbf{U}_{%
\left[ nt\right] },\mathbf{\tilde{V}}_{\left[ nt\right] }(r_{1},r_{2}),S_{n}%
\right) e^{-S_{n}};L_{n}\geq 0\right] }{\mathbf{E}\left[ e^{-S_{n}};L_{n}%
\geq 0\right] } \\
&\rightarrow &\iiiint \Theta ^{\ast }\left( \alpha ,\beta ,\lambda ;u,%
\mathbf{v},-z\right) \mathbf{P}^{+}\left( B^{+}\in du\right) \mathbf{P}%
_{z}^{-}\left( \mathbf{V}_{\infty }(r_{1},r_{2})\in d\mathbf{v}\right) \nu
_{1}\left( dz\right) .
\end{eqnarray*}
\end{lemma}

Now we consider the function 
\begin{equation*}
\omega \left( \lambda ;\mathbf{u},v,x\right) :=\frac{l(v)}{1+\left( \lambda
l(u_{1})+l(u_{2})l(v)e^{x}\right) \kappa (u_{3})}.
\end{equation*}%
Clearly, $\omega \left( \alpha ;\mathbf{u},v,x\right) $ is continuous and
does not exceed 1. For the particular case $\mathbf{u}=\left( 1,1,u\right) $
we use one more notation 
\begin{equation}
\omega ^{\ast }\left( \lambda ;u,v,x\right) :=\frac{l(v)}{1+\left( \lambda
+l(v)e^{x}\right) \kappa (u)}.  \label{explispsi}
\end{equation}%
With the functions above and $\mathbf{v}=\left( v_{1},v_{2}\right) $ we
associate the functions%
\begin{equation*}
\Omega \left( \lambda ;\mathbf{u},\mathbf{v},x\right)
:=\prod_{i=1}^{2}\omega \left( \lambda ;\mathbf{u},v_{i},x\right) ,\qquad
\Omega ^{\ast }\left( \lambda ;u,\mathbf{v},x\right) :=\prod_{i=1}^{2}\omega
^{\ast }\left( \lambda ;u,v_{i},x\right) .
\end{equation*}

Let now $\mathbf{U}_{\left[ nt\right] }=\left( \tilde{U}_{\left[ nt\right]
}^{(1)},\tilde{U}_{\left[ nt\right] }^{(2)},\tilde{U}_{\left[ nt\right]
}^{(3)}\right) $ and $\mathbf{\tilde{V}}_{\left[ nt\right] }(r_{1},r_{2})$
be the same as in Lemma \ref{LLright} above with one exception: one should
take $\alpha =1$ in the definition of the components of $\mathbf{U}_{\left[
nt\right] }$.

By Lemma \ref{Lend} and relations\ (\ref{DDU}) and (\ref{DDV1}) we see that
the following statement is valid.

\begin{lemma}
\label{LLright}Under conditions A1 to A4 
\begin{eqnarray*}
&&\frac{\mathbf{E}\left[ \Omega \left( \lambda ;\mathbf{U}_{\left[ nt\right]
},\mathbf{\tilde{V}}_{\left[ nt\right] }(r_{1},r_{2}),S_{n}\right)
e^{S_{n}};\tau (n)=n\right] }{\mathbf{E}\left[ e^{S_{n}};\tau (n)=n\right] }
\\
&\rightarrow &\iiiint \Omega ^{\ast }\left( \lambda ;u,\mathbf{v},-z\right) 
\mathbf{P}_{z}^{+}\left( B^{+}\in du\right) \mathbf{P}^{-}\left( \mathbf{V}%
_{\infty }(r_{1},r_{2})\in d\mathbf{v}\right) \mu _{1}\left( dz\right) .
\end{eqnarray*}
\end{lemma}

To go further, observe that if the offspring probability functions are
fractional-linear, then by Lemma \ref{LfrLin} for any $0\leq m<n,s\in
\lbrack 0,1],$ and $0\leq r_{2}<r_{1}<1$ 
\begin{eqnarray}
&&f_{0,m}(sf_{m,n}(r_{1}))-f_{0,m}(sf_{m,n}(r_{2}))  \notag \\
&&\qquad\qquad:=sG_{m,n}(s;r_{1})G_{m,n}(s;r_{2})e^{-S_{n}}\frac{r_{1}-r_{2}%
}{\left( 1-r_{1}\right) \left( 1-r_{2}\right) },  \label{Differen}
\end{eqnarray}%
where%
\begin{eqnarray}
G_{m,n}(s;r) &:=&\frac{1-f_{m,n}(r)}{e^{-S_{m}}+\left( 1-sf_{m,n}(r)\right)
b_{m}}  \notag \\
&=&\frac{\left( 1-f_{m,n}(r)\right) e^{S_{m}-S_{n}}\times e^{S_{n}}}{1+\left[
\left( 1-s\right) e^{S_{m}}+s\left( 1-f_{m,n}(r)\right)
e^{S_{m}-S_{n}}\times e^{S_{n}}\right] b_{m}}  \label{G1} \\
&=&\frac{\left( 1-f_{m,n}(r)\right) e^{S_{m}-S_{n}}}{e^{-S_{n}}+\left[
\left( 1-s\right) e^{S_{m}}e^{-S_{n}}+s\left( 1-f_{m,n}(r)\right)
e^{S_{m}-S_{n}}\right] b_{m}}  \label{G2} \\
&\leq &G_{m,n}(1;r)=1-f_{0,n}(r)\leq 1.  \label{G3}
\end{eqnarray}

Introduce the notation 
\begin{equation*}
Y_{t}^{\left( n\right) }:=Z_{\left[ nt\right] }e^{-S_{\left[ nt\right]
}},\,t\in \left( 0,1\right) .
\end{equation*}%
Clearly, for $\lambda \geq 0$%
\begin{equation*}
\mathbf{E}\left[ e^{-\lambda Y_{t}^{\left( n\right) }};T=n+1\right] =\mathbf{%
E}\left[ e^{-\lambda Y_{t}^{\left( n\right) }}\left( f_{\left[ nt\right]
,n+1}^{Z_{\left[ nt\right] }}\left( 0\right) -f_{\left[ nt\right] ,n}^{Z_{%
\left[ nt\right] }}\left( 0\right) \right) \right] =\mathbf{E}[F_{n,t}\left(
\lambda \right) ],
\end{equation*}%
where, in view of (\ref{Differen}) and with $s=s(\lambda ):=\exp \left\{
-\lambda e^{-S_{\left[ nt\right] }}\right\} $ 
\begin{eqnarray*}
F_{n,t}\left( \lambda \right) &:=&f_{0,\left[ nt\right] }\left( sf_{\left[ nt%
\right] ,n}\left( f(0)\right) \right) -f_{0,\left[ nt\right] }\left( sf_{%
\left[ nt\right] ,n}\left( 0\right) \right) \\
&=&sG_{\left[ nt\right] ,n}(s;f(0))G_{\left[ nt\right]
n}(s;0)e^{-S_{n}}X_{f}\left( 1\right) .
\end{eqnarray*}

\begin{lemma}
\bigskip \label{Lfuncjleft}For each fixed $j$%
\begin{equation*}
\mathcal{A}_{j}^{\ast }(\lambda ):=\lim_{n\rightarrow \infty }n^{3/2}\mathbf{%
E}\left[ F_{n,t}\left( \lambda \right) ;\tau (n)=j\right] =K_{1}\mathbf{E}%
\left[ e^{-S_{j}}I\left( \tau (j)=j\right) X_{f}\left( 1\right) A_{j}^{\ast
}(\lambda )\right] ,
\end{equation*}%
where%
\begin{equation*}
A_{j}^{\ast }(\lambda ):=\iiiint \Theta ^{\ast }\left( a_{j},b_{j},\lambda
;u,\mathbf{v},-z\right) \mathbf{P}^{+}\left( B^{+}\in du\right) \mathbf{P}%
_{z}^{-}\left( \mathbf{V}_{\infty }(f(0),0)\in d\mathbf{v}\right) \nu
_{1}\left( dz\right) .
\end{equation*}
\end{lemma}

\textbf{Proof}. Let, as earlier, $\mathcal{F}_{j}^{\ast }$ be the $\sigma $%
-algebra generated by $Q_{1},..,Q_{j}$ and $Q$. \ Denote%
\begin{equation*}
U_{j,n,t}^{(1)}:=\frac{1-\exp \left\{ -\lambda a_{j}e^{-\hat{S}_{\left[ nt%
\right] -j}}\right\} }{\lambda a_{j}e^{-S_{\left[ nt\right] }}}%
,U_{j,n,t}^{(2)}:=\exp \left\{ -\lambda a_{j}e^{-\hat{S}_{\left[ nt\right]
-j}}\right\} ,U_{j,n,t}^{(3)}:=\hat{b}_{\left[ nt\right] -j},
\end{equation*}%
\begin{equation*}
\tilde{V}_{j,n,t}(r):=\left( 1-\hat{f}_{\left[ nt\right] -j,n-j}(r)\right)
e^{\hat{S}_{\left[ nt\right] -j}-\hat{S}_{n-j}}.
\end{equation*}%
By splitting $S_{n}$\ as $\left( S_{n}-S_{j}\right) +S_{j}$ we deduce from (%
\ref{G2}) and (\ref{G3}) that for $s=s(\lambda )=\exp \left\{ -\lambda
a_{j}e^{-\hat{S}_{\left[ nt\right] -j}}\right\} $ and $r\in \lbrack 0,1)$ 
\begin{eqnarray*}
&&G_{nt,n}(s;r)I\left( \tau (n)=j\right) \\
&&\qquad\overset{d}{=}\frac{\tilde{V}_{j,n,t}(r)}{a_{j}e^{-\hat{S}_{n-j}}+%
\left[ \lambda a_{j}U_{j,n,t}^{(1)}e^{-\hat{S}_{n-j}}+U_{j,n,t}^{(2)}\tilde{V%
}_{j,n,t}(r)\right] (b_{j}+a_{j}U_{j,n,t}^{(3)})}\times \\
&&\qquad\qquad\qquad\qquad\qquad\qquad\qquad\qquad\times I\left( \tau
(j)=j\right) I\left( \hat{L}_{n-j}\geq 0\right) \\
&&\qquad=\theta \left( a_{j},b_{j},\lambda ;\mathbf{U}_{j,n,t},\tilde{V}%
_{j,n,t}(r),\hat{S}_{n-j}\right) I\left( \tau (j)=j\right) I\left( \hat{L}%
_{n-j}\geq 0\right) \leq 1.
\end{eqnarray*}%
Hence it follows that%
\begin{equation*}
\mathbf{E}\left[ F_{n,t}\left( \lambda \right) ;\tau (n)=j\right] =\mathbf{E}%
\left[ e^{-S_{j}}I\left( \tau (j)=j\right) X_{f}\left( 1\right)
A_{j,nt,n}^{\ast }\left( \lambda \right) \right] ,
\end{equation*}%
where%
\begin{equation*}
A_{j,nt,n}^{\ast }\left( \lambda \right) :=\mathbf{E}\left[ \Theta \left(
a_{j},b_{j},\lambda ;\mathbf{U}_{j,n,t},\mathbf{\tilde{V}}_{j,n,t}(f(0),0),%
\hat{S}_{n-j}\right) e^{-\hat{S}_{n-j}};\hat{L}_{n-j}\geq 0|\mathcal{F}%
_{j}^{\ast }\right] .
\end{equation*}%
Using now Lemma \ref{Lfuncjleft}, the asymptotic representation (\ref{Star})
and applying the dominated convergence theorem we see that%
\begin{eqnarray*}
&&\lim_{n\rightarrow \infty }n^{3/2}\mathbf{E}[ F_{n,t}\left(
\lambda\right);\tau (n)=j] \\
&&\qquad\qquad=\mathbf{E}\left[ e^{-S_{j}}I\left( \tau (j)=j\right)
X_{f}\left( 1\right) \lim_{n\rightarrow \infty }n^{3/2}A_{j,nt,n}^{\ast
}\left( \lambda \right) \right] \\
&&\qquad\qquad\quad=K_{1}\mathbf{E}\left[ e^{-S_{j}}I\left( \tau
(j)=j\right) X_{f}\left( 1\right) A_{j}^{\ast }(\lambda )\right] =\mathcal{A}%
_{j}^{\ast }(\lambda ),
\end{eqnarray*}%
as desired.

\begin{lemma}
\label{Lfuncright}For any fixed $j$%
\begin{equation*}
\mathcal{B}_{j}^{\ast }(\lambda ):=\lim_{n\rightarrow \infty }n^{3/2}\mathbf{%
E}\left[ F_{n,t}\left( \lambda \right) ;\tau (n)=n-j\right] =K_{2}\mathbf{E}%
\left[ e^{\hat{S}_{j}}I\left( \hat{L}_{j}\geq 0\right) X_{f}\left( 1\right)
B_{j}^{\ast }(\lambda )\right] ,
\end{equation*}%
where%
\begin{equation*}
B_{j}^{\ast }(\lambda ):=\iiint\Omega ^{\ast }\left( \lambda ;u,\mathbf{v}%
,-z\right) \mathbf{P}_{z}^{+}\left( B^{+}\in du\right) \mathbf{P}^{-}\left( 
\mathbf{V}_{\infty }(\hat{f}_{0,j}(f(0)),\hat{f}_{0,j}(0))\in d\mathbf{v}%
\right) \mu _{1}\left( dz\right) .
\end{equation*}
\end{lemma}

\textbf{Proof}. Let $\mathcal{F}_{n-j+1,n+1}$ be the $\sigma $-algebra
generated by $Q_{n-j+1},..,$ $Q_{n+1}$. Denote%
\begin{equation*}
U_{n,t}^{(1)}:=\frac{1-\exp \left\{ -\lambda e^{-S_{nt}}\right\} }{\lambda
e^{-S_{nt}}},U_{n,t}^{(2)}:=\exp \left\{ -\lambda e^{-S_{nt}}\right\}
,U_{n,t}^{(3)}:=b_{nt},
\end{equation*}%
\begin{equation*}
\tilde{V}_{j,n,t}(r):=\left( 1-f_{nt,n-j}(r)\right) e^{S_{nt}-S_{n-j}}.
\end{equation*}%
By splitting $S_{n}$\ as $\left( S_{n}-S_{j}\right) +S_{j}$ and letting $%
s=s(\lambda )=\exp \left\{ -\lambda e^{-S_{nt}}\right\} $ and $r\in \lbrack
0,1)$ we deduce from (\ref{G1}) and (\ref{G3}) that 
\begin{eqnarray*}
&&G_{nt,n}(s;r)I\left( \tau (n)=n-j\right)
=G_{nt,n-j}(s;f_{n-j,n}(r))I\left( \tau (n)=n-j\right) \\
&&\overset{d}{=}\frac{\tilde{V}_{j,n,t}(\hat{f}_{0,j}(r))}{1+\left[ \lambda
U_{n,t}^{(1)}+U_{n,t}^{(2)}\tilde{V}_{j,n,t}(\hat{f}_{0,j}(r))e^{S_{n-j}}%
\right] U_{n,t}^{(3)}}\,I\left( \tau (n-j)=n-j\right) I\left( \hat{L}%
_{j}\geq 0\right) \\
&=&\omega \left( \lambda ;\mathbf{U}_{n,t},\tilde{V}_{j,n,t}(\hat{f}%
_{0,j}(r)),S_{n-j}\right) I\left( \tau (n-j)=n-j\right) I\left( \hat{L}%
_{j}\geq 0\right) \leq 1.
\end{eqnarray*}%
Hence 
\begin{eqnarray*}
\mathbf{E}\left[ F_{n,t}\left( \lambda \right) ;\tau (n)=n-j\right] &=&%
\mathbf{E}[\mathbf{E}\left[ F_{n,t}\left( \lambda \right) I\left( \tau
(n)=n-j\right) |\mathcal{F}_{n-j+1,n+1}\right] ] \\
&=&\mathbf{E}\left[ se^{\hat{S}_{j}}X_{f}\left( 1\right) I\left( \hat{L}%
_{j}\geq 0\right) B_{n,t,j}^{\ast }\left( \lambda \right) \right] ,
\end{eqnarray*}%
where 
\begin{equation*}
B_{n,t,j}^{\ast }\left( \lambda \right) :=\mathbf{E}\left[ \Omega (\lambda ;%
\mathbf{U}_{n,t},\mathbf{\tilde{V}}_{j,n,t}(\hat{f}_{0,j}(f(0)),\hat{f}%
_{0,j}(0)),S_{n-j})e^{S_{n-j}};\tau (n-j)=n-j\right] .
\end{equation*}%
The needed statement follows now from Lemma \ref{LLright} and the dominated
convergence theorem.

The following lemma is crucial for our subsequent arguments.

\begin{lemma}
\label{onedim} For any $t\in \left( 0,1\right) ,$ as $n\rightarrow \infty $%
\begin{equation*}
\mathcal{L}\left( Y_{t}^{\left( n\right) }|T=n+1\right) \rightarrow \mathcal{%
L}\left( W\right)
\end{equation*}%
weakly, where $W$ is an a.s. positive proper random variable.
\end{lemma}

\textbf{Proof. }Relation (\ref{middl}) gives for sufficiently large $n$ and
all $N\geq N\left( \varepsilon \right) $%
\begin{eqnarray}
\mathbf{E}[ F_{n,t}( \lambda );N<\tau(n)
<n-N]\qquad\qquad\qquad\qquad\qquad\qquad  \notag \\
\leq \mathbf{P}( T=n+1;N<\tau ( n) <n-N) \leq \varepsilon n^{-3/2}.
\label{middl2}
\end{eqnarray}%
By Lemmas \ref{Lfuncjleft} and \ref{Lfuncright} for each fixed $N$%
\begin{equation*}
\lim_{n\rightarrow \infty }n^{3/2}\mathbf{E}\left[ F_{n,t}\left( \lambda
\right) ;\tau \left( n\right) \notin \lbrack N+1,n-N]\right]
=\sum_{j=0}^{N}\left( \mathcal{A}_{j}^{\ast }(\lambda )+\mathcal{B}%
_{j}^{\ast }(\lambda )\right)
\end{equation*}%
which, in view of (\ref{AsymExt}), implies%
\begin{equation*}
H^{\ast }(\lambda )=\mathbf{E}e^{-\lambda W}:=\lim_{n\rightarrow \infty }%
\mathbf{E}\left[ e^{-\lambda Y_{t}^{\left( n\right) }}|T=n+1\right] =\frac{1%
}{H(1)}\sum_{j=0}^{\infty }\left( \mathcal{A}_{j}^{\ast }(\lambda )+\mathcal{%
B}_{j}^{\ast }(\lambda )\right) .
\end{equation*}

It follows from the definitions (\ref{explisphi}) and (\ref{explispsi}) that 
\begin{equation*}
\lim_{\lambda \rightarrow \infty }\Theta ^{\ast }\left( \alpha ,\beta
,\lambda ;u,\mathbf{v},x\right) =0\text{ and }\lim_{\lambda \rightarrow
\infty }\Omega ^{\ast }\left( \lambda ;u,\mathbf{v},x\right) =0
\end{equation*}%
leading by the dominated convergence theorem to $\lim_{\lambda \rightarrow
\infty }H^{\ast }(\lambda )=0.$ Thus, the distribution of $W$ has no atom at
zero. On the other hand, again by the dominated convergence theorem, 
\begin{equation*}
\lim_{\lambda \downarrow 0}\mathcal{A}_{j}^{\ast }(\lambda )=K_{1}\mathbf{E}%
\left[ e^{-S_{j}}I\left( \tau (j)=j\right) X_{f}\left( 1\right) A_{j}^{\ast
}(0)\right] ,
\end{equation*}%
where%
\begin{equation*}
A_{j}^{\ast }(0):=\iiiint \Theta ^{\ast }\left( a_{j},b_{j},0;u,\mathbf{v}%
,-z\right) \mathbf{P}^{+}\left( B^{+}\in du\right) \mathbf{P}_{z}^{-}\left( 
\mathbf{V}_{\infty }(f(0),0)\in d\mathbf{v}\right) \nu _{1}\left( dz\right) .
\end{equation*}%
We know by (\ref{DefVinfin}) that 
\begin{equation*}
\mathbf{V}_{\infty }(f(0),0)=\left( \frac{1}{B^{-}+\left( 1-f(0)\right) ^{-1}%
},\frac{1}{B^{-}+1}\right)
\end{equation*}%
which leads to%
\begin{eqnarray*}
A_{j}^{\ast }(0) &=&\iiint \Theta ^{\ast }\left( a_{j},b_{j},0;u,\frac{1}{%
v+\left( 1-f(0)\right) ^{-1}},\frac{1}{v+1},-z\right) \times \\
&&\qquad\qquad\qquad\qquad\times \mathbf{P}^{+}\left( B^{+}\in du\right) 
\mathbf{P}_{z}^{-}\left( B^{-}\in dv\right) \nu _{1}\left( dz\right) .
\end{eqnarray*}%
Recalling (\ref{DefphiDisc}) and (\ref{explisphi}) we see that for any $h>0$%
\begin{eqnarray*}
\phi \left( ha_{j},b_{j},a_{j};u,v,-z\right) &=&\frac{1}{%
a_{j}he^{z}+b_{j}+a_{j}u+a_{j}e^{z}v} \\
&=&\frac{\left( v+h\right) ^{-1}}{a_{j}e^{z}+\left( v+h\right) ^{-1}\left(
b_{j}+a_{j}u\right) } \\
&=&\theta ^{\ast }\left( a_{j},b_{j},0;\left( v+h\right) ^{-1},-z\right) .
\end{eqnarray*}%
Hence, letting $h=1$ and $h=\left( 1-f(0)\right) ^{-1}$ we obtain%
\begin{equation*}
\Theta ^{\ast }\left( a_{j},b_{j},0;u,\frac{1}{v+\left( 1-f(0)\right) ^{-1}},%
\frac{1}{v+1},-z\right) =\Phi \left( a_{j},\frac{a_{j}}{1-f(0)}%
,b_{j},b_{j},a_{j},a_{j};u,v,-z\right) .
\end{equation*}%
This equality, and representations (\ref{Acalif}) and (\ref{Asimp}) show
that $\lim_{\lambda \downarrow 0}\mathcal{A}_{j}^{\ast }(\lambda )=\mathcal{A%
}_{j}(0)$ for any $j=0,1,....$ In a similar way one can prove that \ $%
\lim_{\lambda \downarrow 0}\mathcal{B}_{j}^{\ast }(\lambda )=\mathcal{B}%
_{j}(0)$ for any $j=0,1,....$ This implies $\lim_{\lambda \downarrow
0}H^{\ast }(\lambda )=H(1)$. Thus, the distribution of $W$ has no atom at
infinity.

The lemma is proved.

\textbf{Proof of Theorem \ref{fucttheor}} . For fixed $\delta \in \left(
0,1/2\right) $ we introduce the process with constant paths%
\begin{equation*}
W_{t}^{n}:=e^{-S_{\left[ n\delta \right] }}Z_{\left[ n\delta \right] },t\in %
\left[ \delta ,1-\delta \right] .
\end{equation*}%
By Lemma \ref{onedim} $\mathcal{L}\left( W_{t}^{n},\delta \leq t\leq
1-\delta \right) \rightarrow \mathcal{L}\left( W\right) $ in distribution in
the space $D\left[ \delta ,1-\delta \right] $ with respect to the Skorokhod
topology. Since the limiting process is continuous, we have convergence in
the metric of uniform convergence as well. To prove the theorem it is
sufficient to show that for all $\varepsilon >0$%
\begin{equation}
\lim_{n\rightarrow \infty }\mathbf{P}\left( \sup_{t\in \left[ \delta
,1-\delta \right] }| Y_{t}^{n}-W_{t}^{n}|>\varepsilon \Big|\, T=n+1\right)
=0.  \label{key}
\end{equation}

\bigskip To simplify the subsequent arguments we introduce for $\delta \in
\left( 0,1/2\right) $ and $M>1$ the events 
\begin{equation*}
\mathcal{D}\left( \varepsilon ,\delta \right) :=\left\{ \sup_{t\in \left[
\delta ,1-\delta \right] }\left\vert Y_{t}^{n}-W_{t}^{n}\right\vert
>\varepsilon \right\} ,\quad \mathcal{K}(M,\delta ):=\left\{ Y_{1-\delta
}^{n}\in \lbrack M^{-1},M]\right\}
\end{equation*}%
and denote by $\mathcal{K}^{c}(M,\delta )$ the event complimentary to $%
\mathcal{K}(M,\delta )$. In view of Lemma~\ref{onedim} for any $\gamma >0$
there exists $M$ such that for all $n\geq n(\gamma ,M)$%
\begin{equation*}
\mathbf{P}\left( \mathcal{D}\left( \varepsilon ,\delta \right) \cap \mathcal{%
K}^{c}(M,\delta ),T=n+1\right) \leq \mathbf{P}\left( \mathcal{K}%
^{c}(M,\delta ),T=n+1\right) \leq \gamma n^{-3/2}.
\end{equation*}%
Besides, by the arguments used to demonstrate (\ref{middl}) one can show
that for any $\gamma >0$ there exists $N$ such that for all $n\geq n(\gamma
,N)$%
\begin{eqnarray*}
&&\mathbf{P}\left( \mathcal{D}\left( \varepsilon ,\delta \right) \cap 
\mathcal{K}(M,\delta ),T=n+1;\tau (n)\in \lbrack N+1,n-N-1]\right) \\
&&\qquad \qquad \qquad \leq \mathbf{P}\left( T=n+1;\tau (n)\in \lbrack
N+1,n-N-1]\right) \leq \gamma n^{-3/2}.
\end{eqnarray*}%
Clearly, for any $m<n$%
\begin{eqnarray*}
R(Z_{m},n):= &&f_{m,n+1}^{Z_{m}}(0)-f_{m,n}^{Z_{m}}(0)\leq
Z_{m}f_{m,n+1}^{Z_{m}-1}(0)\left( f_{m,n+1}(0)-f_{m,n}(0)\right) \\
&& \\
&=&Z_{m}f_{m,n+1}^{Z_{m}-1}(0)\left( 1-f_{m,n+1}(0)\right) \left(
1-f_{m,n}(0)\right) e^{S_{m}-S_{n}}X_{f_{n+1}}\left( 1\right) .
\end{eqnarray*}%
Hence, by the inequalities $1-f_{m,n+1}(0)\leq 1-f_{m,n}(0)\leq
e^{S_{n}-S_{m}}$ and (\ref{EstX}) and Assumption A1 we see that%
\begin{equation}
R(Z_{m},n)\leq \frac{Z_{m}}{e^{S_{m}}}e^{S_{n}}X_{f_{n+1}}\left( 1\right)
\leq \rho \frac{Z_{m}}{e^{S_{m}}}\,e^{S_{n}}.  \label{RR1}
\end{equation}%
On the other hand, by the inequality $1-x\leq e^{-x},\,x>0,$ we get%
\begin{equation*}
f_{m,n+1}^{Z_{m}-1}(0)=\left( 1-\left( 1-f_{m,n+1}(0)\right) \right) ^{Z_{m}}%
\frac{1}{f_{m,n+1}(0)}\leq e^{-Z_{m}\left( 1-f_{m,n+1}(0)\right) }\,\frac{1}{%
f_{m,n+1}(0)}
\end{equation*}%
which, in view of the estimates 
\begin{equation}
f_{m,n+1}(0)\geq f_{m+1}(0)\geq \chi ,\quad \frac{1-f_{m,n}(0)}{%
1-f_{m,n+1}(0)}\leq \frac{1}{1-f_{m+1}(0)}\leq \frac{1}{1-\chi }\leq 2,\quad
\label{Cr1}
\end{equation}%
gives%
\begin{eqnarray*}
&&R(Z_{m},n) \\
&\leq &\frac{Z_{m}}{e^{S_{m}}}e^{2S_{m}}(1-f_{m,n+1}(0))^{2}e^{-Z_{m}\left(
1-f_{m,n+1}(0)\right) }\frac{e^{-S_{n}}}{f_{m,n+1}(0)}\frac{1-f_{m,n}(0)}{%
1-f_{m,n+1}(0)}X_{f_{n+1}}\left( 1\right) \\
&&\qquad \leq 2\rho \chi ^{-1}\frac{Z_{m}}{e^{S_{m}}}\,e^{-S_{n}}\sup_{x\geq
0}x^{2}\exp \left\{ -\frac{Z_{m}}{e^{S_{m}}}x\right\} .
\end{eqnarray*}%
For $m=\left[ n(1-\delta )\right] $ inequalities (\ref{RR1}) and (\ref{Cr1})
give 
\begin{equation}
R(Z_{\left[ n(1-\delta )\right] },n)I\left( \mathcal{K}(M,\delta )\right)
\leq \rho Me^{S_{n}}  \label{Remright}
\end{equation}%
and 
\begin{eqnarray}
R(Z_{\left[ n(1-\delta )\right] },n)I\left( \mathcal{K}(M,\delta )\right)
&\leq &2\rho \chi ^{-1}Me^{-S_{n}}\sup_{x\geq 0}x^{2}e^{-xM^{-1}}  \notag \\
&\leq &8\rho \chi ^{-1}M^{3}e^{-2}e^{-S_{n}}.  \label{Remleft}
\end{eqnarray}%
Recalling that by Lemma 3.8 in \cite{ABKV}%
\begin{equation}
\mathbf{P}\left( \mathcal{D}\left( \varepsilon ,\delta \right) |\Pi \right)
\leq \left( \varepsilon ^{-2}\left[ \sum_{i=n\delta }^{\left[ n\left(
1-\delta \right) \right] }\eta _{i+1}e^{-S_{i}}+e^{-S_{\left[ n\left(
1-\delta \right) \right] }}-e^{-S_{\left[ n\delta \right] }}\right] \right)
\wedge 1=:\mathcal{U}_{n}  \label{3.21}
\end{equation}%
we get by means of (\ref{Remleft}) for a fixed $j\in \mathbb{N}_{0}$ 
\begin{eqnarray*}
\Xi _{(j)}(n):= &&\mathbf{P}\left( \mathcal{D}\left( \varepsilon ,\delta
\right) \cap \mathcal{K}(M,\delta ),T=n+1;\tau (n)=j\right) \\
&=&\mathbf{E}\left[ I\left( \mathcal{D}\left( \varepsilon ,\delta \right)
\cap \mathcal{K}(M,\delta )\right) \left( f_{\left[ n(1-\delta )\right]
,n+1}^{Z_{\left[ n(1-\delta )\right] }}(0)-f_{\left[ n(1-\delta )\right]
,n}^{Z_{\left[ n(1-\delta )\right] }}(0)\right) ;\tau (n)=j\right] \\
&\leq &8\rho \chi ^{-1}M^{3}e^{-2}\mathbf{E}\left[ I\left( \mathcal{D}\left(
\varepsilon ,\delta \right) \right) e^{-S_{n}};\tau (n)=j\right] \\
&\leq &8\rho \chi ^{-1}M^{3}e^{-2}\mathbf{E}\left[ \mathcal{U}%
_{n}e^{-S_{n}};\tau (n)=j\right] .
\end{eqnarray*}

\bigskip By Lemma 3.1 in \cite{ABKV}, for all $x\geq 0$ 
\begin{equation*}
\mathcal{U}_{n}\rightarrow 0\;\;\;\;\;\;\;\;\text{\ \ }\mathbf{P}_{x}^{+}-%
\text{a.s.}
\end{equation*}%
Hence, applying the arguments similar to those used in the proof of Lemma %
\ref{Lfuncjleft} we see that, as $n\rightarrow \infty ,$ $\ \Xi
_{(j)}(n)=o(n^{-3/2})$ for each fixed $j$. Further, using inequality (\ref%
{Remright}) we get%
\begin{eqnarray*}
\Xi ^{(j)}(n):= &&\mathbf{P}\left( \mathcal{D}\left( \varepsilon ,\delta
\right) \cap \mathcal{K}(M,\delta ),T=n+1;\tau (n)=n-j\right) \\
&&\qquad \qquad \leq \rho M\mathbf{E}\left[ \mathcal{U}_{n}e^{S_{n}};\tau
(n)=n-j\right] .
\end{eqnarray*}%
Applying now arguments similar to those used to demonstrate Lemma \ref%
{Lfuncright} one can show that as $n\rightarrow \infty ,$ $\Xi
^{(j)}(n)=o(n^{-3/2})$ for each fixed $j$.

\bigskip\ Combining the estimates above and taking first the limit as $%
n\rightarrow \infty $ then as $\gamma \downarrow 0$ and, finally, as $%
M\rightarrow \infty $ we arrive at (\ref{key}).

The theorem is proved.

\end{document}